\documentclass[12pt]{article}
\usepackage{amsmath,amssymb,amsthm,amscd,graphicx}
\usepackage{latexsym}
\usepackage{enumerate}
\usepackage{amsbsy, amsfonts, textcomp}
\usepackage{makeidx}
\newtheorem{theorem}{Theorem}[]
\newtheorem{proposition}[theorem]{Proposition}

\theoremstyle{definition}
\newtheorem{definition}[theorem]{Definition}
\newtheorem{example}[theorem]{Example}

\theoremstyle{remark}
\newtheorem{remark}[theorem]{Remark}
\newcommand{\Aut}{\mbox{Aut }}

\newcommand{\G}{\mathbb{G}}
\newcommand{\R}{\mathbb{R}}
\newcommand{\C}{\mathbb{C}}

\newcommand{\N}{\mathbb{N}}
\newcommand{\cL}{\mathcal{L}}

\def \DGal{\operatorname{DGal}}
\def \DAut{\operatorname{DAut}}
\def \Aut{\operatorname{Aut}}
\def \GL{\operatorname{GL}}
\def \SO{\operatorname{SO}}
\begin{document}

\Large
\centerline{\bf Real Liouville extensions}

\large
\centerline{Teresa Crespo, Zbigniew Hajto}

\let\thefootnote\relax\footnotetext{Both authors acknowledge support of grant MTM2009-07024, Spanish Science Ministry.}

\normalsize

\vspace{0.5cm}
\begin{abstract}
We give a characterization of real Liouville extensions by differential Galois groups.
\end{abstract}

\section{Introduction}

In \cite{CHS1} and \cite{CHS2}, we proved the existence of a Picard-Vessiot extension for linear differential equations defined over a real differential field with a real closed field of constants $C$, we gave an appropriate definition of its differential Galois group, proved that it has the structure of a $C$-defined linear algebraic group and established a Galois correspondence theorem in this setting.
In \cite{GK}, Gel'fond and Khovanskii characterized Liouville functions over $\R$ using  differential rings of real functions with a finiteness property.

In this paper, we give a positive answer to the following question posed by Askold Khovanskii: ``Is it true that a necessary and sufficient condition for
solvability of a real differential equation by real Liouville functions
follows from real Picard-Vessiot theory?'' More precisely, we characterize linear differential equations defined over a real differential field with a real closed field of constants $C$, which are solvable by real Liouville functions, as those having a differential Galois group whose identity component is solvable and $C$-split.

All fields considered will be of characteristic zero. We refer the reader to \cite{CH} for the topics on differential Galois theory, to \cite{BCR} for those on real fields.

\section{Preliminaries}

In this section, we recall the definitions and main properties of differential field extensions and linear algebraic groups used in the sequel.

\begin{definition} Let $L|K$ be a differential field extension, $\alpha$ an element in $L$. We say that $\alpha$ is

\begin{enumerate}[-]
\item \emph{an integral over $K$} if  $\alpha'= a \in K$ and $a$ is not a derivative in $K$;
\item \emph{the exponential of an integral over $K$} if $\alpha'/\alpha \in K\setminus \{0\}$.
\end{enumerate}
\end{definition}

\begin{definition} A differential field extension $K\subset L$ is called a \emph{Liouville extension} (resp. a \emph{generalised Liouville extension}) if there exists a chain of intermediate differential fields $K=F_1 \subset F_2 \subset \dots \subset F_n =L$ such that $F_{i+1}=F_i(\alpha_i)$, where $\alpha_i$ is either an integral or the exponential of an integral over $F_i$ (resp. or $\alpha_i$ is algebraic over $F_i$).
\end{definition}

\begin{definition} Let $K$ be a real differential field with real closed field of constants $C$, $\cL(Y)=0$ a homogeneous linear differential equation defined over $K$, $y$ a solution to $\cL(Y)=0$. We shall say that $y$ is a \emph{real Liouville solution to} $\cL(Y)=0$ if $y$ is contained in some real Liouville extension of $K$.
\end{definition}

\begin{definition} (\cite{CHS2} Section 4.1)\label{dg} Let $K$ be a real differential field with real closed field of constants $C$, $L|K$ a  real differential field extension. The differential Galois group $\DGal(L|K)$ of $L|K$ is the set $DHom_K(L,L(i))$ of $K$-~differential morphisms from $L$ into $L(i)$ endowed with the group structure translated from the group $\DAut_{K(i)} (L(i))$ of $K(i)$-differential automorphisms of $L(i)$ by means of the bijection

$$\begin{array}{ccc} \DAut_{K(i)} (L(i)) & \rightarrow & DHom_K(L,L(i)) \\ \tau & \mapsto & \tau_{|L} \end{array}.$$
\end{definition}

\begin{proposition} (\cite{CHS2} Prop. 4.1) Let $K$ be a real differential field with real closed field of constants $C$, $L|K$ a  real Picard-Vessiot extension. The differential Galois group $\DGal(L|K)$ of $L|K$ has the structure of a $C$-defined linear algebraic group.
\end{proposition}

\begin{remark} Let $\overline{C}$ denote the algebraic closure of the field $C$. An affine variety $V$ in the affine space $\overline{C}^n$ is $C$-defined if and only if $V^{\sigma}=V$, for all $\sigma \in \Aut_{C}(\overline{C})$. From this fact, we deduce easily that if $G$ is a $C$-defined linear algebraic group, then its identity component $G^0$ is $C$-defined as well.

\end{remark}

We study now the extensions of a real differential field with real closed field of constants obtained by adjunction of an integral and those obtained by adjunction of the exponential of an integral.

\begin{example}\label{integ}
Let $K$ be a real differential field with real closed field of constants $C$. We consider the differential extension $L=K\langle \alpha
\rangle$, obtained by adjunction of an integral. As in the case in which $C$ is algebraically closed, one may prove
that $\alpha$ is transcendental over~$K$, $K\langle \alpha
\rangle|K$ is a Picard-Vessiot extension, and $\DGal(K\langle \alpha
\rangle|K)$ is isomorphic to the additive group, which is $C$-defined. Since $K\langle \alpha
\rangle|K$ is a transcendental extension, the field $K\langle \alpha
\rangle$ may be ordered.
\end{example}

\begin{example}\label{prim} Let $K$ be a real differential field with real closed field of constants $C$. We consider the differential extension $L=K\langle \alpha
\rangle$, obtained from $K$ by adjunction of the
exponential of an integral. We assume that the field $L$ is a real field and that its field of constants is equal to $C$. Then $K\langle \alpha \rangle|K$ is a Picard-Vessiot extension. As in the case in which $C$ is algebraically closed, one may prove
that if $\alpha$ is algebraic over $K$, then $\alpha^n
\in K$ for some $n \in \N$. The Galois group $\DGal(L|K)$ is
isomorphic to the multiplicative group if $\alpha$ is
transcendental over $K$ and to a finite cyclic group if $\alpha$ is
algebraic over $K$. These groups are $C$-defined.
\end{example}

The next proposition gives the classification of the connected linear algebraic groups of dimension~1 defined over a real closed field $C$. Besides the additive and multiplicative groups appearing in the algebraically closed case, we also have the circle, which is a $C$-form of the multiplicative group.

\begin{proposition}
Let $G$ be a connected linear algebraic group of dimension~1 defined over a real closed field $C$. Then $G$ is isomorphic over $C$ either to $\G_m$, $\G_a$ or $\SO_2$.
\end{proposition}

\begin{example}\label{cos} Let us consider the differential equation $Y''+Y=0$ defined over the field $\C$ of complex numbers. Its solutions lie in the field $\C(e^{it})$, which is a Liouville extension of $\C$. However, if we consider the same equation as defined over $\R$ its differential Galois group is $\SO_2$ and it will follow from theorem \ref{th} that it does not have real Liouville solutions.
\end{example}

\begin{definition} Let $G$ be a connected solvable linear algebraic group defined over a field $C$. We say that $G$ is $C$-split if it has a composition series $G=G_1 \supset G_2 \supset \cdots \supset G_s = 1$ consisting of connected $C$-defined closed subgroups such that $G_i/G_{i+1}$ is $C$-isomorphic to $\G_a$ or $\G_m$, $1\leq i <s$.
\end{definition}

\begin{proposition}\label{LK} (\cite{B}, 15.4)  Let $G$ be a connected linear algebraic group of degree $n$, defined over the field $C$. The following conditions are equivalent:
\begin{enumerate}
\item $G$ is solvable and $C$-split.
\item $G$ is solvable and all rational characters of $G$ are defined over $C$.
\item there exists an element $x \in \GL(n,C)$ such that $xGx^{-1} \subset Tr_n$.
\end{enumerate}
\end{proposition}

\section{Main results}

In this section we shall characterize linear differential equations defined over a real differential field with real closed field of constants admitting a fundamental set of real Liouville solutions.

\begin{proposition}\label{propo} Let $K$ be a real differential field with real closed field of constants $C$. Let $L$ be a real Picard-Vessiot extension of $K$, which is a Liouville extension of $K$. Then the differential Galois group $\DGal(L|K)$ is solvable and $C$-split.
\end{proposition}

\noindent {\it Proof.} Clearly, $L(i)$ is a Liouville extension of $K(i)$ and the algebraically closed field $C(i)$ is the field of constants of both $K(i)$ and $L(i)$. By definition \ref{dg} and \cite{CH}, Prop. 6.4.2, it only remains to prove that $G=\DGal(L|K)$ is $C$-split. By hypothesis, $L|K$ has a chain of intermediate differential fields $K=F_1 \subset F_2 \subset \cdots \subset F_s =L$ such that $F_{i+1}$ is obtained from $F_i$ either by adjunction of an integral or by adjunction of the exponential of an integral. By the fundamental theorem of real Picard-Vessiot extensions, $G$ has a normal chain of closed subgroups $G \supset \DGal(L|F_2) \supset \cdots \supset 1$ and $\DGal(L|F_i)/\DGal(L|F_{i+1}) \simeq \DGal(F_{i+1}|F_i)$ is $C$-isomorphic to $\G_m$ or $\G_a$, $1\leq i <s$.  \hfill $\Box$

\begin{definition} A real differential field extension is \emph{r-normal} if for any $x \in L\setminus K$, there exists $\sigma \in \DGal(L|K)$ such that $\sigma(x) \neq x$.
\end{definition}

\begin{remark} A real Picard-Vessiot extension is r-normal.
\end{remark}

\begin{proposition}\label{tri}
Let $L|K$ be an r-normal extension of real
differential fields. Assume that there exist elements
$u_1,\dots,u_n \in L$ such that for every $\sigma \in \DGal(L|K)$ we have

\begin{equation}\label{eq:liou}
\sigma \, u_j = a_{1j} \, u_1 + \dots + a_{j-1,j} \, u_{j-1}+
a_{jj} \, u_j \, , \, j=1,\dots,n,
\end{equation}

\noindent with $a_{ij}$ constants in $L(i)$ (depending on $\sigma$).
Then $K\langle u_1,\dots,u_n \rangle$ is a Liouville extension of
$K$.

\end{proposition}

\noindent {\it Proof.} The proof follows the same steps as for \cite{CH} Proposition 6.4.3.

\begin{theorem} Let $K$ be a real differential field with real
closed field of constants $C$. Let $L$ be a real Picard-Vessiot
extension of $K$. Assume that the identity component $G^0$ of
$G=\DGal(L|K)$ is solvable and $C$-split. Then $L$ can be obtained from $K$ by a
finite r-normal extension followed by a Liouville extension.
\end{theorem}

\noindent {\it Proof.} Let $F=L^{G^0}$.
Then $F|K$ is a finite extension, $F(i)|K(i)$ is a normal extension and $\DGal(L|F) \simeq
G^0$. Then by Proposition \ref{LK}, we can apply Proposition
\ref{tri} and obtain that $L|F$ is a Liouville extension.

\begin{theorem}\label{th}
 Let $K$ be a real differential field with real
closed field of constants $C$. Let $L$ be a Picard-Vessiot
extension of $K$. Assume that $L$ can be embedded in a real
differential field $M$ which is a generalized Liouville extension
of $K$ with no new constants. Then the identity component $G^0$ of
$G=\DGal(L|K)$ is solvable and $C$-split.
\end{theorem}

\noindent {\it Proof.} Clearly $L(i)$ is a Picard-Vessiot extension of $K(i)$ and $M(i)$ is a generalized Liouville extension of $L(i)$. By the characterization of Liouville extensions of differential fields with algebraically closed field of constants (see e.g. \cite{CH} Theorem 6.5.4), $G^0$ is solvable and $L(i)$ may be obtained from $K(i)$ by a finite normal extension, say $F|K(i)$, followed by a Liouville extension, which is a Picard-Vessiot extension with differential Galois group $G^0$. Hence $L$ is obtained from $K$ by the finite r-normal extension $F\cap L|K$ followed by the extension $L|F\cap L$, which is a real Liouville extension. By Proposition \ref{propo}, $G^0$ is $C$-split.

\begin{example} Let us consider the differential equation $Y''+Y=0$ defined over the field $\R$ of real numbers. As stated in Example \ref{cos}, its differential Galois group is $\SO_2$. The map

$$\begin{array}{ccc} \SO_2 & \rightarrow & \G_m \\ \left( \begin{array}{cc} a & -b \\ b & a \end{array}\right) & \mapsto & a + b \, i \end{array}
$$

\noindent is clearly a rational character of $\SO_2$ which is not $\R$-defined, hence the equation $Y''+Y=0$ has no real Liouville solutions.
\end{example}

\vspace{1cm}
\footnotesize
\begin{tabular}{lcl} 
Teresa Crespo && Zbigniew Hajto \\
Departament d'\`{A}lgebra i Geometria && Faculty of Mathematics and Computer Science \\  Universitat de Barcelona && Jagiellonian University \\ Gran Via de les Corts Catalanes 585 && ul. Prof. S. \L ojasiewicza 6 \\ 
08007 Barcelona, Spain && 30-348 Krak\'ow, Poland \\
teresa.crespo@ub.edu && zbigniew.hajto@uj.edu.pl
\end{tabular}
\end{document}